\documentclass[11pt]{amsart}
\usepackage{amscd}
\newtheorem{thm}{Theorem}
\newtheorem{prop}{Proposition}
\newtheorem{lem}{Lemma}

\theoremstyle{remark}

\theoremstyle{definition}
\newtheorem{defn}{Definition}

\newcommand{\N}{\mathbb{N}}

\newcommand{\F}{\mathbb{ F}}

\newcommand{\Z}{\mathbb{ Z}}

\newcommand{\Hom}{\operatorname{Hom}}

\newcommand{\lra}{\longrightarrow}

\newcommand{\Q}{\mathbb{ Q}}
\newcommand{\R}{\mathbb{ R}}

\newcommand{\im}{{\rm im\,}}

\newcommand{\ra}{{\rightarrow}}
\newcommand{\Zp}{\mathbb{Z}/p}
\newcommand{\Fp}{\mathbb{F}_p}

\title{Poincar\'e duality in P.A. Smith Theory}

\author{Christopher~Allday}
\address{C.~Allday, University of Hawaii, Department of Mathematics, 2565
Mc Carthy Mall, Honolulu, HI}
\email{chris@math.hawaii.edu}

\author{Bernhard~Hanke}
\address{B.~Hanke, Universit\"at M\"unchen, Theresienstr. 39, 80333
M\"unchen, Germany}
\email{hanke@rz.mathematik.uni-muenchen.de}

\author{Volker~Puppe}
\address{V.~Puppe, Universit\"at Konstanz, 78457 Konstanz, Germany}
\email{Volker.Puppe@uni-konstanz.de}

\begin{document}

\begin{abstract} Let $G=S^1$, $G=\Z / p$ or more generally $G$ be a finite
$p$-group, where $p$ is an odd prime. If $G$ acts on a space  whose
cohomology ring
fulfills Poincar\'e duality (with appropriate coefficients $k$), we prove a
mod $4$ congruence
between  the total Betti number of $X^G$ and a number which depends only on
the $k[G]$-module
structure of $H^*(X;k)$. This improves the well known mod $2$ congruences
that hold for
actions on general spaces.
\end{abstract}

\date{\today; MSC 2000: Primary 57S10, 57P10, 55N10; Secondary 55N91}
\thanks{The second author is holding a DFG research grant. He would like to  thank the University of Notre Dame for its hospitality during
the work on this paper.}

\maketitle

\section{Introduction and statement of results}

Let $X$ be a finite dimensional connected CW complex and let $k$ be a
commutative ring with unit. We
say that $X$ is a {\em Poincar\'e duality space} over $k$ (we will
often write {\em $k$-PD space} instead) of {\em
formal dimension} $n$, if
$H^*(X;k)$ is a finitely generated $k$-module and if there is a class
$\nu\in H_n(X;k)$
such that
\[
  H^*(X;k) \ra H_{n-*}(X;k), \quad c \ra c \cap \nu
\]
is an isomorphism. Note that if $k$ is a field, this is equivalent to
requiring
that
\[
   H^*(X;k) \times H^*(X;k) \stackrel{\cup}{\lra} H^*(X;k)
\stackrel{\nu}{\lra} k
\]
is a nonsingular pairing (viewing $\nu$ as an element in $\Hom
(H^{*}(X;k),k)$).

Let $G$ denote the group $S^1$ (with its usual topology) or
$\Z/p$, where $p$ is an  odd prime number. Let $k=\Q$, if $G=S^1$, and
$k=\Fp$, if $k=\Z/p$. By a well
known result proven independently by Chang-Skjelbred in
\cite{CS} and Bredon in \cite{Br}, each component of the fixed point set of
a finite
dimensional $G$-CW complex $X$ fulfills Poincar\'e duality over $k$, if
this is property holds for
$X$. By now there are several further versions and variants of proofs of
this result (\cite{AP, Ha1,
Ra}). In this paper,  we will use certain consequences of Poincar\'e
duality to deduce relations between
the total Betti number of $X^G$ and the $k[G]$-module $H^*(X;k)$.

\begin{thm} \label{th1} Let $G = S^1$ and let $X$ be a finite dimensional
connected $G$-CW complex
such that $X$ is a $\Q$-PD space of formal dimension $n$. If
\begin{itemize}
  \item $n$ is even or
  \item $n=2m+1$, $X^G \neq \emptyset$, $H^i(X;\Q) = 0$ for $ 0 < i \leq
m$, $i$ even,
\end{itemize}
the following congruence holds.
\[
     \dim_{\Q} H^*(X^G;\Q) \equiv \dim_{\Q} H^*(X;\Q) \mod 4\, .
\]
\end{thm}

If $G = S^1 \times \ldots \times S^1$ and if $X$ is a finite  dimensional
connected
$G$-CW complex  with finitely many orbit types that fulfills Poincar\'e
duality over $\Q$, then an
analogue of Theorem
\ref{th1} holds. This is true, because in this case we can choose a
subcircle $S^1 \subset G$, such
that $X^{S^1} = X^G$.

A version of Theorem \ref{th1} for actions of $G= \Z/p$, where $p$ is an
odd prime, can be proven under
additional assumptions on the  space $X$. It turns out that one also has to
take into account the
fact that now the induced
$G$-action on $H^*(X;\Fp)$ might be nontrivial. Using a spectral
sequence argument together with results from \cite{Ha2}, we show:

\begin{thm} \label{th2} Let $G =\Zp$ (where $p$ is an odd prime) and let
$X$ be a finite dimensional
connected
$G$-CW complex such that $X$ is an $\Fp$-PD space of formal dimension $n$.
Furthermore, assume that $H^*(X;\Z_{(p)})$ does not contain $\Zp$ as 
a direct summand.  Then,  we get a decomposition as graded $\Fp[G]$-modules
\[
    H^*(X;\Fp) = F^* \oplus T^* \oplus R^* \, ,
\]
where $F^*$ is a free $\Fp[G]$-module, $T^*$ is a trivial $\Fp[G]$-module and
$R^*$ is a direct sum of $\Fp[G]$-modules of the form $\ker \epsilon$, where
$\epsilon:\Fp[G]\ra \Fp$ is the augmentation map. If
\begin{itemize}
 \item $n$ is even or
 \item $n = 2m+1$, $X^G \neq \emptyset$, $T^i = 0$ for $0 < i \leq m$, $i$
even, $R^i = 0$ for $0 < i
\leq m$, $i$ odd,
\end{itemize}
the following congruence holds.
\[
   \dim_{\Fp} H^*(X^G;\Fp) \equiv \dim_{\Fp} T^* + \frac{1}{p-1} \dim_{\Fp}
R^* \mod 4\, .
\]
\end{thm}

Theorem \ref{th1} was first proved by A.~Sikora in his PhD thesis. This and Theorem \ref{th2} for the 
case that $R^* = 0$ and $H^*(X;\Z)$ does not contain any $p$-torsion are also contained in his paper
\cite{AS}, where a somewhat different line of argument is used.

It is well known that 
the cited result by Bredon, Chang and Skjelbred
immediately generalizes
to actions of finite $p$-groups $G$. This is achieved by applying
induction on the
order of $G$ and using the fact that every finite (nontrivial) $p$-group
contains a normal subgroup of
order $p$. This procedure can be applied in our situation as well
and shows the following result.

\begin{thm} \label{th3} Let $G$ be a finite $p$-group (where $p$ is
an odd prime) and let $X$ be a finite dimensional connected $G$-CW
complex that is an $\Fp$-PD space of even
formal dimension. If $p > \dim H^*(X;\Fp)$, then
\[
 \dim_{\Fp} H^*(X^G;\Fp) \equiv \dim_{\Fp} H^*(X;\Fp) \mod 4 \, .
\]
\end{thm}

If we
impose more restrictions on
$X$ and on the group operation, then  using a result from \cite{Ha3}, we
can use a
similar induction argument for comparing rational Betti numbers:

\begin{thm} \label{th4} Let $G$ be a finite $p$-group (where $p$ is an odd
prime) and let $X$ be
a finite connected simplicial complex with $G$ acting simplicially such
that $X$ is an even
dimensional orientable $\Z_{(p)}$-homology manifold. Furthermore, assume
that $p > \dim
H^*(X;\Fp)$. Then
\[
  \dim_{\Q} H^*(X^G;\Q) \equiv \dim_{\Q} H^*(X;\Q)   \mod 4 \, .
\]
\end{thm}

In the case of smooth actions, this result can be proven without making use
of \cite{Ha3}.

In Section 5, we will show by examples that none of the additional
assumptions in the second part of Theorem \ref{th1} and Theorem \ref{th2}
can be dropped.

We would like to thank the referee
for some helpful comments and suggestions.

\section{Algebraic preliminaries}

For our later applications to spectral sequences, it will be useful to provide a notion of Poincar\'e duality for bigraded algebras.

\begin{defn} \label{poinc} Let $k$ be a field and let $A^{*,*}$ be a
$(\Z/2\times \N)$-bigraded
associative and graded commutative $k$-algebra with unit
that  is finitely generated as $k$-vector space.  (The 
graded commutativity is required with respect to the total $\Z
/2$-grading, where $A^{i,j}$ has total degree $i + j \mod 2$.) 
\begin{itemize}
\item[i.] We call $A^{*,*}$ {\em connected}, if $A^{0,0} \cong k$.
\item[ii.] $A^{*,*}$ is a {\em $k$-Poincar\'e duality algebra} of {\em
formal dimension} $n$,
if there is a surjective linear map  $\phi:A^{*,*} \ra k$ (called {\em
orientation}) with $\phi(A^{i,j}) = 0$ for $(i,j) \neq (0,n)$ 
such that the bilinear form
\[
  \zeta: A^{*,*} \times A^{*,*} \stackrel{{\rm mult.}}{\lra} A^{*,*}
\stackrel{\phi}{\lra} k
\]
is nondegenerate.
\end{itemize}
\end{defn}

Note that in a connected Poincar\'e duality algebra of formal dimension
$n$, we have $A^{0,n} \cong k$.
The following elementary fact provides a connection between the  Euler
characteristic and  the total
Betti number of Poincar\'e duality algebras of even formal dimension.

\begin{lem} \label{l1} Let $A^{*,*}$ be a Poincar\'e duality algebra of even
formal dimension over a field $k$. If ${\rm char~} k\neq 2$, then
\[
   \dim_k A^{*,*} \equiv \chi(A^{*,*})  \mod 4 \, ,
\]
where  the Euler characteristic is calculated using the total
$\Z/2$-grading of $A$.

\end{lem}

\begin{proof} Let $2m$ be the formal dimension of $A^{*,*}$. Using the
 induced total $\Z/2$ grading of $A^{*,*}$, we claim
that the dimension of
$A^{odd}$ is even, which obviously proves  the lemma. For all $i$, we get
isomorphisms
\begin{eqnarray*}
  A^{ 0, 2i+1} & \cong & \Hom (A^{0,2m-2i-1},k) \\
  A^{1, 2i}    & \cong & \Hom (A^{1,2m-2i},k)
\end{eqnarray*}

by Poincar\'e duality. Hence,  $A^{0,2i+1}$ and $A^{0, 2m-2i-1}$, respectively
$A^{1,2i}$ and $A^{1,2m-2i}$ have the same
dimension. Furthermore, if $m$ is even, the module $A^{1,m}$, carries a skew
nonsingular form by Poincar\'e duality and therefore has even dimension, as
${\rm char\,} k\neq 2$.
The same is true  for
$A^{0,m}$, if $m$ is odd.
\end{proof}

\begin{prop} \label{p1} Let $(A^{*,*},\delta)$ be a connected $(\Z/2\times
\N)$-graded differential
algebra over the field $k$ with a differential $\delta$ of arbitrary bidegree,
but lowering the second grading parameter  and acting as a
derivation. Furthermore, assume that
$A^{*,*}$ is a Poincar\'e duality algebra of formal dimension
$n$. Then, after taking homology, we get either
$H(A^{*,*},\delta)=0$ or $H(A^{*,*},\delta)$ is again a connected
Poincar\'e duality algebra of formal
dimension $n$.
\end{prop}

\begin{proof} Assume that $H(A^{*,*}) \neq 0$. This implies that the unit
of $A$ (that
generates  $A^{0,0}$) is not in the image of $\delta$. 
Let $c \in A^{0,n} \cong k$ be a generator. Assume that
$\delta(c) = x \neq 0$.  By Poincar\'e duality, we then find an element
$y \in A^{*,*}$  with $c=  yx \in A^{0,n}$. We therefore get
\[
  0 \neq x = \delta(c)  = \delta(yx) = \delta (y)  x \, ,
\]
which implies that $\delta(y)$ is a generator of $A^{0,0}$, contrary to
what we said before.
Hence, because $c$ is not hit by $\delta$, any orientation $\phi$ of $A$ induces a surjective linear map
$H^{*,*}(A) \ra k$. To complete the
proof, observe that by the derivation property of $\delta$, we have
\[
   \zeta(\ker \delta, \im \delta) =0\, ,
\]
using the bilinear form $\zeta$ from Definition \ref{poinc}.  As
\[
  \dim_k \ker \delta + \dim_k \im \delta = \dim_k A^{*,*}\, ,
\]
we may conclude that $\ker \delta$ is exactly the orthogonal complement of
${\rm im~} \delta$
with respect to $\zeta$. This proves that the induced bilinear form on
$H^{*,*}(A)$ is
nonsingular.
\end{proof}

\begin{prop} \label{p2} Let $(A^{*,*},\delta)$ be as in the last
Proposition. Additionally, we
assume that ${\rm char~}k \neq 2$,  the formal dimension of $A^{*,*}$ is an
odd number $2m+1$, the differential $\delta$ has odd (total) degree,
lowering the second grading
paramter in $A^{*,*}$ and  $A^{0,i} = 0$ for $0 < i \leq m$, $i$ even,
$A^{1,i} =0$ for $0 < i \leq m$, $i$ odd. If $H(A^{*,*}) \neq 0$. Then
\[
   \dim_{k} A^{*,*} \equiv \dim_{k} H(A^{*,*},\delta)  \mod 4\, .
\]
\end{prop}

\begin{proof} Let $Z^*= \ker \delta \subset A^*$ denote the cycles in
$A^*$ with respect to $\delta$ (here we use the total $\Z/2$ grading,
again). By our assumption on
$A^*$,  the even dimensional part $H^{ev}(A^*,\delta)$ of the homology  of
$A^*$ with respect to
$\delta$, coincides with $Z^{ev}$, cf.~the proof of Proposition \ref{p1}.
Now we define a  bilinear form
\[
 \gamma: A^{ev} \times A^{ev} \ra k, \quad (x,y) \ra \phi( x \cdot
\delta(y))\, ,
\]
using an orientation $\phi$ of $A$. It is  easy to check, that $\gamma$
induces a well defined,
nonsingular and skew symmetric form on
$A^{ev}/Z^{ev}$, hence this vector space has even dimension over $k$ (using
the fact that ${\rm char~}k \neq 2$). As
the Euler characteristic of
$A^*$ and of
$H(A^*)$ are equal, we get the equation
\[
   A^{ev} - H^{ev}(A^*) = A^{odd} - H^{odd}(A^*)\, .
\]
Therefore, the number $\dim A^* - \dim H^*(A^*)$ is divisible by $4$.
\end{proof}

\section{Proof of Theorems \ref{th1}, \ref{th3} and \ref{th4}}

Here, when applying
the results from the last section, we usually forget about the first grading
parameter of $A^{*,*}$ and use the  $\N$-grading by the second parameter.
Part of the
argument  is based on the following well known fact (cf. \cite{AP},
Exercise (3.29)):

\begin{lem} \label{l2} Let $G=S^1$ or $G= \Z/p$ and let $X$ be a finite
dimensional $G$-CW complex such
that
$H^*(X;\Z)$  is a finitely generated $\Z$-module. Then
\begin{itemize}
   \item[i.]   $\chi(X^G) = \chi(X)$, if $G=S^1$.
   \item[ii.]   $\chi(X^G) = \Lambda(g)$, if $G=\Z/p$ and $g \in G$, $g
\neq 1$.
\end{itemize}
Here, $\Lambda(g)$ denotes the Lefschetz number
\[
    \Lambda(g) = \sum_i (-1)^i {\rm trace} (g_* \colon H_i(X;\Q) \ra
H_i(X;\Q)) \, ,
\]
regarding $g$ as a map $X \ra X$.
\end{lem}

If we set $A^{0,1} = 0$ and  $A^{*,0} = H^*(X;\Q)$, 
respectively $A^{*,0} = H^*(X^G;\Q)$,   then Lemma \ref{l1} and
Lemma \ref{l2} give  the following sequence of 
equations:
\[
   \dim H^*(X^G;\Q) \equiv \chi(X^G) = \chi(X) \equiv \dim H^*(X;\Q)
{\rm~~mod~} 4\, .
\]
This shows the first part of Theorem \ref{th1}.

The proof of Theorem \ref{th3} proceeds by induction on the order of
$G$. Assume $|G|\neq
1$ and choose a normal subgroup $H \subset G$, $H\cong \Zp$ that exists by
group theory.
Each component of $X^H$ is an $\Fp$-PD space of even formal dimension by
the Theorem of Bredon-Chang-Skjelbred.
As $\dim H^*(X;\Q) < p$ by our assumption on $p$, the induced action
of $H$ on $H^*(X;\Q)$ is trivial: Let $V$ be a rational vector
space of dimension smaller than $p-1$ and let $f$ be a linear 
endomorphism of $V$ with $f^p = {\rm id}$. Then the minimal polynomial
of $f$ in $\Q[x]$ must divide $x^p -1$. The last polynomial splits over $\Q[x]$ into $(x-1)$ and an irreducible factor of degree $p-1$. Because the minimal polynomial of $f$ has degree at most $\dim V$, 
it must therefore be equal to $x-1$. This covers all cases, where $H^*(X;\Q)$ is not concentrated
in degree $0$. In the remaining case, the assertion follows from the fact that $H^0(X;\Q)$ is a permutation module. 

Altogether, for $1 \neq h \in H$, we obtain
\[
    \Lambda(h) = \chi(X) \, .
\]

By Lemma \ref{l1} and
Lemma \ref{l2} above, we get

\begin{eqnarray*}
   \dim H^*(X^H;\Fp)  \equiv   \chi(X^H) =
   \chi(X)  \equiv \dim H^*(X;\Fp) {\rm~~mod~} 4\, .
\end{eqnarray*}

(Note that the Euler characteristic does not depend
on the coefficient field used).
For the induction step, observe that by Smith theory
\[
  \dim H^*(F;\Fp) \leq \dim H^*(X;\Fp) < p
\]
for each component $F$ of $X^G$. Furthermore, the group $G' = G/H$ has
order less than the order
of $G$ and each component of $X^H$ is invariant under the induced
$G'$-action by our assumption
on $p$. Using the Theorem of  Bredon-Chang-Skjelbred, each
component of $X^H$ is again an $\Fp$-PD space of even
formal dimension. Hence the induction
hypothesis applies
to each component of $X^H$. The proof of Theorem \ref{th3} is now
complete.

For the proof of Theorem \ref{th4}, we recall the following fact.

\begin{prop} (\cite{Ha3}, Proposition 13) Let $G = \Zp$ act simplicially
on a finite simplicial
complex $X$ that is an orientable $\Z_{(p)}$-homology manifold. If $F
\subset X$ is a component of the fixed point set $X^G$, then $F$
is an orientable $\Z_{(p)}$-homology manifold of even
codimension in $X$.
\end{prop}

Using this fact, Theorem \ref{th4} follows from Theorem \ref{th3} by
observing that the total
Betti numbers with either $\Fp$ or $\Q$ coefficients of a space are
congruent modulo $4$, if
this space fulfills
Poincar\'e duality  both over $\Q$ and over $\Fp$ and has even  formal
dimension for
both fields of coefficients.  This follows by the independence of the Euler
characteristic of the coefficient field and by using Lemma \ref{l1} twice.

\medskip

For the  proof of the second part of Theorem \ref{th1} we consider the
cohomological Leray-Serre spectral
sequence $(E_r,\delta_r)$ for the Borel fibration
\[
     X \ra X_G= EG \times_G X \ra BG
\]
with coefficients in $\Q$. For the $E_2$-term of this spectral sequence as
a module
over $H^*(BG;\Q) \cong \Q[t]$, where $t\in H^2(BG;\Q)$ is a generator, we get
\[
    E_2^{* ,\mu} \cong H^{\mu}(X;\Q) \otimes \Q[t] \, .
\]
Because all the differentials in the spectral sequence are $\Q[t]$-linear,
we can evaluate the terms
$E_r$,  $r\geq 2$ at $t=1$ and because evaluation at $t=1$ commutes with
taking homology with respect
to each $\delta_r$ (see \cite{AP}, Lemma (A.7.2)), we get a
converging spectral sequence $(\overline{E}_r, \overline{\delta}_r)$, where
\[
   \overline{E}_r= (E_r)_{t=1}
\]
and $\overline{\delta}_r$ is induced by $\delta_r$. This spectral sequence
is concentrated in
the first column and $\overline{E}_{\infty}^{*,*}$ is (noncanonically and
not preserving the grading)
isomorphic to
$H^*(X_G;\Q)_{t=1}$ as a $\Q$-vector space. By
the evaluation theorem (cf.~\cite{AP}, Theorem (3.5.1)), we have a
canonical isomorphism (of ungraded modules)
\[
   H(X_G;\Q)_{t=1} \cong H(X^G;\Q) \, .
\]
Because $X$ fulfills Poincar\'e duality over $\Q$, the term
$\overline{E}_2^{*,*}=\overline{E}_2^{0,*}$
is an
$\N$-graded $\Q$-Poincar\'e duality algebra of formal dimension $2m+1$, where
we use the grading induced by $E_2$. As $X^G \neq \emptyset$, we have
$\overline{E}_r \neq 0$ for
all $r$.
Now the claim follows from Propositions  \ref{l1} and \ref{l2} in the second
section of this paper.

\section{Proof of Theorem \ref{th2}}

Let $p$ be an odd prime number, let $G =\Zp$  and let  $X$ be a  finite
dimensional connected $G$-CW
complex that  is a $\Fp$-PD space of formal dimension $n$. Furthermore, we
assume that
$\beta(H^*(X;\Fp)) = 0$, where
\[
   \beta\colon H^*(X;\Fp) \ra H^{*+1}(X;\Fp)
\]
is the Bockstein operator associated to the exact sequence of coefficients
\[
  0 \ra \Z/p \ra \Z / p^2 \ra \Z/p \ra 0 \, .
\]
This condition is equivalent to the requirement, that each cohomology class
in $H^*(X;\Fp)$ can
be lifted to a cohomology class in $H^*(X;\Z /p^2)$, which is the case, if
and only if $H^*(X;\Z_{(p)})$ does not contain a direct summand of the form $\Z /p$.
Note that the $G$-action
on $X$  induces a $\Z /p^2 [G]$-module structrue on $H^*(X;\Z /p^2)$. By
our assumption on $\beta$,
$H^*(X;\Z /p^2)$ is a free $\Z /p^2$-module, so the following proposition
is an immediate consequence of \cite{Ha2}, Proposition 6, which is  
shown using a little representation theory over $\Z /p^2$.

\begin{prop} \label{decomp} For each $i$ there is a $\Fp[G]$-linear
decomposition
\[
    H^i(X;\Fp) = H^i(X; \Z /p^2) \otimes \Fp \cong T^i \oplus F^i \oplus
R^i \, ,
\]
where $T^i$ is a trivial $\Fp[G]$-module, $F^i$ is a direct sum of free
$\Fp[G]$-modules
and $R^i$ is a direct sum of modules of the form $\ker \epsilon$,  where
$\epsilon:\Fp[G] \ra \Fp$
is the augmentation map.
\end{prop}

Let $(E_r,\delta_r)$ denote the cohomological Leray-Serre spectral sequence
for the Borel fibration
\[
     X \ra X_G= EG \times_{G} X \ra BG
\]
with coefficients in $\Fp$. For $r \geq 2$, each $E_r$ is a differential
bigraded algebra over
$H^*(BG;\Fp) \cong \Fp[t] \otimes \Lambda(s)$, where
$t \in H^2(BG;\Fp)$ and $s \in H^1(BG;\Fp)$ are generators and $\beta(s)
=t$. By the evaluation theorem
(cf.~\cite{AP}, Theorem (1.4.5)), the inclusion $X^G \hookrightarrow X$
induces an isomorphism of (ungraded)
$\Lambda(s)$-algebras
\[
   H(X_G;\Fp)_{t=1} \cong H(X^G;\F_p) \otimes \Lambda(s)\, ,
\]
hence, after evaluating at $s=0$, we get an induced isomorphism
\[
   H(X_G;\Fp)_{t=1,s=0} \cong H(X^G;\F_p) \, .
\]
Now we encounter the following difficulty that did not arise in the
consideration of $S^1$-actions
above: In
order to calculate the dimension of the left hand side of the second
isomorphism one has to get around
the difficulty that evaluation at
$s=0$ does not  commute with taking homology in general, and so a spectral
sequence
argument as in  Section 3 does not seem to be feasible in this case.
However, under
the additional  assumption $\beta=0$, we can get around this difficulty by
using the following fact (cf.~\cite{Ha2}, Proposition 9, and 
the last part of the proof of Proposition 10).

\begin{prop} \label{free} For all $r\geq 2$, the
localized terms
$E_r[t^{-1}]$ are  finitely generated free $(\Z \times \N)$-bigraded
differential $\Lambda(s) \otimes
\Fp[t,t^{-1}]$-algebras.  Furthermore, evaluation  at $t=1$ and $s=0$  on $E_r[t^{-1}]$
commutes with taking homology with respect to the differentials induced
by $\delta_r$.
\end{prop}

Now, we set
\[
  \overline{E}_r = (E_r)_{t=1,s=0}
\]
and use the induced differentials $\overline{\delta}_r$. Each term
$\overline{E}_r$  has an induced
$(\Z/2 \times \N)$-grading and for $r=2$ we obtain
\[
    \overline{E}_2^{\gamma,\mu} \cong
H^{\gamma}(G;H^{\mu}(X;\Fp))_{t=1,s=0} \cong \left\{
\begin{array}{ll}
          \Fp^{\dim T^{\mu}}, {\rm~~if~} \gamma {\rm~~is~even}, \\
          \Fp^{\frac{\dim R^{\mu}}{p-1}}, {\rm~~if~} \gamma {\rm~~is~odd}\, .
         \end{array} \right.
\]
Here, we used Proposition \ref{decomp}
and the fact (which follows from usual dimension shifting) that
$H^i(\Fp;\ker \epsilon)[t^{-1}]$ is isomorphic to $\Fp$ in odd
degrees and is equal to zero in
even degrees. Further, we get from Proposition \ref{free} and the
evaluation theorem
\[
  \dim_{\Fp} \overline{E}_{\infty} = \frac{1}{2}
\dim_{\Fp}(E_{\infty})_{t=1} =
  \dim_{\Fp} H^*(X^G;\Fp) \, .
\]
Now, the proof of Theorem \ref{th3} is completed by using induction in the
spectral sequence.
Notice that $\overline{E}_2^{*,*}$ is a $(\Z/2\times \N)$-bigraded
connected Poincar\'e duality algebra
over $\Fp$ of formal dimension $n$ in the  sense of Definition \ref{poinc}
(cf.~\cite{Ha2}, proof of Proposition 9). For the induction step, if
$n$ is even, we use Lemma
\ref{l1} and the fact that the  Euler characteristic of a $\Z/2$-graded
differential complex does not
change after  taking homology with respect to a differential of odd degree.
If $n$ is odd, we use
Proposition \ref{p1} and Proposition \ref{p2}.

\section{Examples, applications and concluding remarks}

The following examples illustrate the significance of the conditions 
on the Betti numbers in the second part of Theorem \ref{th1} and
\ref{th2}. The example of free $S^1$ or $\Zp$-actions
on spheres of odd dimension shows that the requirement $X^G\neq \emptyset$
in the second part of Theorem
\ref{th1}
and in the second part of Theorem \ref{th2} is necessary.
In \cite{Br2}, p.~425, an example of an $S^1$-action on $X=S^3 \times S^5
\times S^9$ is constructed whose fixed point set is an $S^7$-bundle over
$S^3 \times
S^5$ and has total
Betti number $6$ (with coefficients  in $\Q$). As the total Betti number of
$X$ is $8$,
one sees that even if fixed points exist, the
assumption on
the vanishing of certain Betti numbers in the second part of Theorem \ref{th1}
is needed. Restricting this $S^1$-action to $\Z/p\subset S^1$, we similarly
may conclude
that all additional assumptions in the second part of Theorem \ref{th2} are
necessary.

Now, we will construct examples of Poincar\'e duality spaces of odd formal
dimension that
fulfill the additional requirement on Betti numbers in the second part
of Theorem \ref{th1} and \ref{th2}. Let $X$ be an arbitrary connected finite
simplicial
complex with the property that its even dimensional integral cohomology is
concentrated in degree $0$. Now embed $X$ in a Euclidean space $\R^{2m+1}$,
where $2m+1  \geq 2 \dim X+1 $, and take a regular neighbourhood $R$ of $X$ inside
$\R^{2m+1}$ which
can be assumed to be a compact oriented smooth manifold with boundary.
Gluing two
copies
of this manifold (the orientation of one of which had been reversed) along
their boundaries
yields a $(2m+1)$-dimensional connected oriented smooth manifold $Y$, whose
even dimensional integral
cohomology below  dimension $m+1$ is concentrated in degree $0$. 
This follows from the Mayer-Vietoris sequence and the 
fact that by general position, the inclusion $\partial R \hookrightarrow R$ 
induces isomorphisms of homotopy groups up to degree 
$2m+1 - \dim X -2 \geq m - 1$ and a surjection in degree
$2m+1 - \dim X - 1 \geq m$. In particular, the 
connecting homomorphism $H^i(\partial R;\Z)\ra H^{i+1}(Y;\Z)$ 
in the Mayer-Vietoris sequence is $0$, if $i \leq m - 1$. Hence, for
actions on $Y$, the
second part of Theorem \ref{th1} and \ref{th2} can be applied (assuming
that the
induced action on $H^*(Y;\Fp)$ is trivial in the case of $G =\Zp$).

Another example
illustrating the second part of Theorem \ref{th1} and Theorem \ref{th2} can
be constructed as follows. Consider the space $X = S^1 \times S^{2m}$
equipped with
the $S^1$-action that acts trivially on the first factor and is the usual
rotation action (fixing north
and south pole) on the second factor. The fixed point set of this
$S^1$-action is the union of
two circles. Next we choose one point in each fixed point component and
remove small $S^1$-invariant
neighborhoods equivariantly diffeomorphic to $D^{2m+1}$ (with a linear
$S^1$-action) around each of these
two fixed points which gives an
$S^1$-manifold $Z$ with two boundary components each of which is
equivariantly diffeomorphic to $S^{2m}$
with  the rotation action by $S^1$.  Now we form the equivariant connected
sum of $Z$ and $[0,1] \times
S^{2m}$,  where on the last space $S^1$ acts trivially on the first factor
and by a rotation
action on the second. In this way, we obtain a $(2m+1)$-dimensional  oriented
$S^1$-manifold which
fulfills the requirement of the second part of Theorem \ref{th1} and
\ref{th2} (using the induced action by  $\Z /p \subset S^1$). It is easy to
check
that the integral cohomology of this space has rank $1$ in degrees $0$ and
$2m+1$,
rank $2$ in degrees $1$ and $2m$ and is $0$ in all other degrees. The
total Betti number of the fixed
set (which is just a
single copy of
$S^1$) is two. So the  Leray-Serre spectral sequence for the Borel
construction  does not
collapse at the $E_2$-level in this case.

In this paper we have been working in the category of $G$-CW complexes. But
using
$\check{\rm C}$ech cohomology and the usual somewhat more technical
machinery (see,
e.g.~\cite{AP}), one
can extend all the results to general $G$-spaces which fulfill the
hypothesis (LT)
for the localization theorem (see \cite{AP}, p.~208). As this generalization
is straightforward, we leave it to the interested reader.

\end{document}